%% file: pldi16.tex
\documentclass[pldi,english]{article}

\usepackage[T1]{fontenc}
\usepackage{babel}
\usepackage{SIunits}            
\usepackage{courier}            
\usepackage[scaled]{helvet} 
\usepackage{url}                  
\usepackage{listings}          
\usepackage{enumitem}      
\usepackage[colorlinks=true,allcolors=blue,breaklinks,draft=false]{hyperref}   

\begin{document}

%
%

\title{Empirical Research and Automatic Processing Method of Precision-specific Operation}

\author{
  Wang, Ran\\
  \texttt{lilianwangran@163.com} 
  \and
  He, Xinrui\\
  \texttt{hexinrui@pku.edu.cn}
}

\maketitle
\begin{abstract}
  Because of the digit limitation of floating point, 
  significant inaccuracy often occurs during the process of mathematical calculation, 
  which may lead to catastrophic loss. 
  Normally, people believe that adjustment of floating-point precision 
  is an effective way to solve this problem, 
  since high-precision floating-point number has more digits to store information. 
  Thus, it is a prevalent method to reduce the inaccuracy in much floating-point related research, 
  that performing all the operations with higher precision. 
  However, we discover that some operations may lead to larger error in higher precision. 
  In this paper, we define this kind of operation 
  that generates large error due to precision adjustment a precision-specific operation. 
  Furthermore, we propose a light-weight searching algorithm for 
  detecting precision-specific operations and figure out an automatic processing method to fixing them.
  In addition, we conducted an experiment on the scientific mathematic library of GLIBC. 
  The result shows that there are many precision-specific operations, 
  and our fixing approach can significantly reduce the inaccuracy.
\end{abstract}

\input{Introduction}

\input{Background}

\input{Approach}
\input{Implementation}

\input{Research}

\input{Conclusion}

\section{Acknowledgements}

This document is based heavily on ones prepared for previous
conferences and we thank their program chairs; in particular, Sandhya
Dwarkadas (ASPLOS'15), Sarita Adve (ASPLOS'14), Steve Keckler
(ISCA'14), Christos Kozyrakis (Micro'13), Margaret Martonosi
(ISCA'13), Onur Mutlu (Micro'12), and Michael L. Scott (ASPLOS'12).

\bibliographystyle{abbrvnat}


\end{document}

%% file: Introduction.tex
\section{Introduction}

Regardless of whether or not aware of that the widely use of floating-point number operations, 
floating-point numbers have penetrated into almost all the area of computer science. 
Floating-point numbers, numbers that include a decimal fraction, 
are always used to approximate real numbers in mathematics. In most programming language, 
floating-point number format is designed as a necessary element because of 
the indispensable status of floating-point numbers in scientific calculation. 
From personal computer to super computer, from C to Java, 
we could see floating-point operations everywhere.

Unfortunately, floating-point numbers could not represent the exact value of some real numbers 
because of the limitation of the storage or digit format of computer. 
Although the inaccuracy of a single number may be slight, 
the accumulation of error after a series of floating-point numbers operations, 
may be huge enough to create unexpected consequence. 
For instance, a program, aiming at calculate a space shuttle’s orbit, 
involves millions of floating-point number calculations. 
If the programmer does not concern about the inaccuracy of floating-point number operations, 
the accumulated error may lead the aircraft off course or even crashed.

\paragraph{Related work}

To resolve the aforementioned problem, 
recent research works employ the techniques of precision adjusting. 

Increasing precision is a common approach to ensure the calculation 
results more close to the “true” real number in the scientific computation. 
For instance, climate modeling requires complicated mathematical computation and simulation. 
Priori work \cite{Yun:2000} employed a double-double arithmetic for a climate modeling code, 
which reduce the numerical variability of the entire system dramatically. 
By using higher precision, they significantly improve numerical reproducibility 
and stability in Parallel Applications. 
Another research is in the supernova simulation area. 
P. H. Hauschildt and E. Baron \cite{Hauschildt:99} used double-double (128-bit or 31-digit) 
and quad-double (256-bit or 62-digit) types solve the expanding stellar atmosphere problem successfully.
Extending the length of digits to store a more exact number, in another word, 
processing the program in higher precision, is a popular approach to decrease the potential error.

Sometimes, speed is more important than accuracy to an application. 
In this case, researchers prefer to decreasing precision to perform more operations in one second. \cite{Rubio:12} \cite{Schkufza:14}
For instance, \cite{Lam:13} proposed a framework that performs better by decreasing precision. 
They designed a breadth first search algorithm to automate identification of code regions 
that can use lower precision. 
After that, they used binary instrumentation and modification to build mixed-precision configuration 
and get some transformations of the original program. 
Then, they kept decreasing the precision of different transformations until the result is qualified enough to pass the efficiency test. 
Their work limits the accuracy loss of decreasing precision.

\paragraph{Precision-specific operation}
The techniques of precision adjusting base on two assumptions:
\begin{itemize}[noitemsep]
\item Increasing precision leads to more accurate result
\item Decreasing precision usually leads to small accuracy loss, which is acceptable, and higher efficiency.
\end{itemize}

Notwithstanding, the assumptions are invalid in some situation. 
There is a kind of operation could only be calculated in a specific precision level 
and the adjustment of precision could cause huge error. 
Such as an operation is accurate in low precision but lose its accuracy in high precision. 
We define this kind of operation that generates large error 
due to precision adjustment a precision-specific operation. 
We discovered two typical patterns of precision-specific operation:
particular constant and Union.

\paragraph{Contributions}
\begin{itemize}[noitemsep]
\item We propose a methodology of detecting precision-specific computation. 
\item We realize a methodology of fixing precision-specific computation problem.
\item We perform an empirical study of precision-specific computation in standard C library and GSL
\item Our evaluation shows that our method could improve the calculation accuracy significantly.
\end{itemize}

%% file: Background.tex
\section{Background}

\subsection{IEEE Standard for Floating-Point Arithmetic Standard}

In the IEEE Standard for Floating-Point Arithmetic Standard (IEEE 754 Standard), 
the floating-point number has a format below:
 $(- 1)^{s} \times c \times b^{q}$
This format include three integers: s = a sign (zero or one), 
c = a significand (or 'coefficient'), q = an exponent. 
b = 2 for floating point. 
In addition, IEEE Standard for Floating-Point Arithmetic defines 4 types of precision: 
single, single-extended, double and double-extended. 
Table~\ref{table:components} shows the layout for single (32-bit) and double (64-bit) precision floating-point values. 
The number of bits for each field are shown and bit ranges are in square brackets.

\begin{table}[t!]
  \centering
{  \sffamily\small 
  \begin{tabular}{ccccc}
    \textbf{Precision} & \textbf{Sign} & \textbf{Exponent} & \textbf{Fraction} & \textbf{Bias}\\
    \hline
    \textbf{Single} & 1 [31] & 8 [30-23] & 23 [22-00] & 127\\
    \textbf{Double} & 1 [63] & 11 [62-52] & 52 [51-00] & 1023\\
  \end{tabular}
}
  \caption{Single and double precision components. }
  \label{table:components}
\end{table}

\subsection{Error of floating-point number}

The errors sources can be separated into three groups: rounding, data uncertainty, and truncation. \cite{Benz:12} 
Rounding errors are unavoidable due to the finite precision. 
Data uncertainty comes from the initial input or the result of a previous computation. 
Because of the limitation of technology, instrument or approach, 
input data from measurements or estimations is usually inaccurate to some extent. 
Truncation arises when a numerical algorithm approximates a mathematical function. 
Such as a simulation algorithm of Taylor Expansion.
Furthermore, due to the approximation, floating-point arithmetic is prone to accuracy problems 
caused by accumulation and catastrophic cancellation.

\paragraph{Accumulation}

Accumulation could amplify the calculation error. 
Because of the finite precision, some real number could only be approximated, such as real number $0.1$. 
In single precision, the number is storage as $0X3DCCCCCD$, which is the nearest $0.1$ number. 
However, in decimal, it equals to 
$0.10000000149011611938…$
If you add 0.1 for 10000 times,
the result in practice would be 

$999.90289306640625000000000000000$

rather than $1000$. 
The accumulation enlarges the error.

\paragraph{Catastrophic cancellation}

When you try to operate cancellation on two similar numbers, 
you may get catastrophic error of the result. 
For example, you expect variable a equals to $1$ in real number computation, 
however, because of the former accumulated error, 
the value of a is $1.0004$ practically in single precision. 
Then you execute the operation: 
$result = a - 1$
The value of $result$ turns to $0.0004$ rather than 0. 
Here, the relative error is huge.

\subsection{FPdebug}
We use FPdebug \cite{Benz:12} as a tool to analysis floating-point accuracy in program dynamically. 
Fpdebug uses binary translation to perform every floating-point computation 
side by side in higher precision. 
Furthermore, it uses a lightweight slicing approach to track the evolution of errors.

%% file: Approach.tex
\section{Approach}

\subsection{Assumption}

We assume the major difference between precision-specific operation and other operations 
is that precision-specific operation could always produce large error regardless of the value of input,
while other operations seldom or only produce large error to specific input. 
This assumption is reasonable because the error of precision-specific operation 
comes from the discrepancy between the precision the developer set and the actual precision. 
At any value of input, the large error occurs with high probability. 
In particular, the accumulation error is not in our consideration. 
In the aforementioned example, huge error occurs after $10000$ times of addition. 
Nevertheless, there is no significant error in a single addition. 
Thus this addition operation is not precision-specific operation.
Based on this assumption, we propose a light-weight searching algorithm 
for detecting precision-specific operations and an approach for fixing this problem.

\subsection{Definitions}

In this paper, we use relative error to represent the calculation error 
for floating-point number computation. 
The value of a floating-point before the adjustment in precision is called 
original value, or approximate value. 
The value after the adjustment in precision is called shadow value. 
The value of mathematical real number is called exact value. 
We have the formula of relative error below:

$Relative\;error = | \frac{exact\;value - approximate\;value}{exact\;value} |$

Thus, in the aforementioned example which add $0.1$ for $10000$ times, 
the original value of result equals to 

$999.90289306640625000000000000000$. 

The shadow value of result is $1000$. 
Then we get the relative error equals to 
$9.710693359375 \times 10^{-5}$.

\subsection{Approach for detecting precision-specific computation}

In this section, we propose a light-weight method for detecting precision-specific computation. 
Firstly, we generate large amount of input values in the domain of definition 
for each function to be detected. 
Secondly, we transform each operation of the function into three-address instruction format. 
Thirdly, we monitor the relative error of all the floating-point left values in assignment statements. 
With the vast input values, we could get plentiful relative errors. 
Finally, we can locate the first stable operation with large relative error through statistical analysis. 
It is the precision-specific computation we found. 
Here, the “stable” means the operation produces relative errors over a threshold 
for a certain percentage of the input values. 
For instance, the function is executed $n$ times while this operation is executed $m$ times 
($n < m$, $n > m$, $n = m$ are all possible). 
We set a threshold for relative error to $e0$ and a threshold for percentage to $p0$. 
The relative error of left value is greater than $e0$ for $k$ times. 
If $k / m > p0$, this operation could be a precision-specific operation. 
Since $k / m$ is meaningless when m is too small, 
we exclude the situation when $m / n$ less than a threshold $p1$.

\subsection{Approach for fixing precision-specific computation}

In this section, we describe an original method of fixing precision-specific computation. 
Considering the definition of precision-specific computation, 
the inaccuracy happens because of the adjustment of the precision. 
Thus, one approach to solve this problem is to reuse the original precision 
for the precision-specific computation and 
change back to the adjusted precision after the computation is executed. 
Our approach focuses on the precision adjustment on a single instruction 
and fixes the precision-specific computation after detection to some extent. 
Above all, we implement two simple functions to 
support the fixing of precision-specific computation problem:
\begin{itemize}[noitemsep]
\item The reducePrec function

The first step is to test if the shadow value exists. 
If it does, we fetch the original precision level of the variable. 
The second step is to extend or round the shadow value to the original precision 
and assign to the original value.
\item The resumePrec function

This function is an opposite process to the reducePrec function. 
The first step is to test if the shadow value exists. 
If it does, we fetch the original precision level of the variable. 
The second step is to extend or round the original value to the adjusted precision 
and assign to the shadow value.
\end{itemize}
Hence, our method solves the problem in the following process.
After we detect the instruction is a precision-specific computation, 
we add the reducePrec function right before the instruction 
and add the resumePrec function right after the instruction. 
In detail, we mark all the precision-specific computation in the detection process 
and execute the modified program in two precision levels simultaneously.
When the program reaches the marked point, 
instead of executing the precision-specific computation directly, 
it will call the resumePrec function, which transfers the shadow value to the original value. 
Then the precision-specific computation would be calculated in the original precision. 
After the instruction is executed, the program will call the reducePrec function, 
which transfers the original value to the shadow value.

%% file: Implementation.tex
\section{Implementation of approaches}

\subsection{Generation of input values}

Generate at least $1000$ input values in the domain of definition. 
For example, we generate $1000$ input values with range $[-1, 1]$ 
with interval of $0.002$ for function acos.

\subsection{Adjustment of precision}

In practice, we implement the method with Valgrind and Fpdebug. 
Valgrind is an efficient framework for dynamic analysis 
and we use Valgrind to track the calculation results of each instruction. 
Fpdebug \cite{Benz:12} is a tool based on Valgrind, which support multi-precision debug. 
In addition, Fpdebug calculates with MPFR library, 
a C library for multiple-precision floating-point computations with correct rounding. 
We complement Fpdebug in purpose of recording original value, 
shadow value and realizing reducePrec function and resumePrec function, which mentioned before.

\subsection{Transformation of source code}

We use LLVM and CLANG to transform the source code into three-address instruction format. 
This technology is based on analysis of Abstract Syntax Tree (AST). 
The complex operation could be separated into a sequence of three-address instructions 
by implementing the post-order traversal of the Abstract Syntax Tree. 
Besides, in order to monitor the relative error of left value for each assignment statement, 
we insert three functions: reducePrec, resumePrec and computErr for each operation. 
Function computer is for computation of relative error. 
Function reducePrec is for precision decreasing while function resumePrec is for precision increasing. 
Both support float, double and long double as input type. 

In a normal process, we only insert computer for the transformed code. 
When we find a precision-specific operation by detecting approach, 
we insert reducePrec before it and resumePrec after it to fix the problem. 
Then we keep recompiling the modified code until all the precision-specific operations are fixed. 
We insert these three functions to all operations 
because we want to skip the excessively long recompiling time. 
In addition, we introduce switch mechanism to control the execution of three inserted functions. 
Every group of the three functions has a unique switch number. 
For example, if $x = y + z$ has switch number $1000$, 
we get the transformation with switch mechanism below:

$reducePrec(\&x,\;1000);$

$x\; =\; y\; +\; z;$

$resumePrec(\&x, 1000);$

$computeErr(“x”, \&x, 1000);$

We stored the switch numbers of open status in a hash table. 
Each function look up the hash table to decide whether to be executed with time complexity $O(1)$.

%% file: Research.tex
\section{Empirical Research}

\subsection{Set up the research}

\subsubsection{Set up the standard of computation results}
We set up the standard calculation result 
with a reliable multi-precision calculation library created by Professor 
\href{http://www.zhaoshizhong.org/download.htm}{Shizhong Zhao}
in order to evaluate the result of our empirical study of precision specific computation 
in scientific mathematic library of GLIBC. 
This library supports reliable multi-precision operation including addition, subtraction, 
multiplication, division, trigonometric functions, 
inverse trigonometric functions, exponential function and logarithmic function. 
The reliability of this calculation library has been testified by many examples 
which have errors in professional calculators, such as Casio, Microsoft and Google calculator, 
or professional calculation software, such as Matlab. 
For instance, when we calculate $20^{65} – e^{65 \times ln(20)}$, 
the standard result is $0$. 
However, we get the following results in Table~\ref{table:calculators}.

\begin{table}[t!]
  \centering
{  \sffamily\small 
  \begin{tabular}{rl}
    \textbf{Professional Calculator} & \textbf{Calculation Result} \\
    \hline
    Casio (fx-82ES PLUS) & $2.79 \times 10^{72}$\\
    Microsoft (Windows 7) & $-7.30057727803363947538724610 \times 10^{48}$\\
    Google calculator& $1.0007532 \times 10^{71}$\\
    Matlab (R2012a) & $6.1253 \times 10^{70}$\\
    Reliable calculation library & $0$ \\
  \end{tabular}
}
  \caption{Calculation results for common professional calculators. }
  \label{table:calculators}
\end{table}

The reliable calculation library returns $0$, which is correct. 

Here is another example: when we calculate the value of $exp(-0.0277)$, 
the manual calculation result with Taylor expansion for $99$ items in double precision is

$0.972680127073139888516095652449$ 

and we get 

$9.72680127073140 * 10^{-1}$ with GSL double precision. 
However, if we extend the precision level for every instruction inside the function $exp$, 
the calculation result is 

$9.72681569708922 * 10^{-1}$, which has a larger error than GSL double precision. 
The reliable calculation library returns 

$0.972680127073139846902979085281$, when keeps 30 decimal places. 
The reliable calculation library performs more accurate in all the examples we test. 
Hence, we complement mathematical formulas for GLIBC mathematical functions 
according to the GLIBC Manual with the help of reliable calculation library. 
The formulas are listed below in Table~\ref{table:formulas}.

\begin{table}[t!]
  \centering
{  \sffamily\small 
  \begin{tabular}{cc}
    \textbf{Function} & \textbf{Formula} \\
    \hline
    $acos$ & $arccos(a)$ \\
    $acosh$ & $ln(a + (a^{2} - 1)^{\frac{1}{2}})$\\
    $asin$ & $arcsin(a)$ \\
    $asinh$ & $ln(a + (a^{2} + 1)^{\frac{1}{2}})$\\
    $atan$ & $arctan(a)$\\
    $atan2$ & $arctan(\frac{a}{c})$ or $arctan(\frac{a}{c}) + \pi$ \\
          & or $arctan(\frac{a}{c}) - \pi$ \\
    $atanh$ & $\frac{1}{2} \ times ln(\frac{(1+a)}{(1-a)})$\\
    $cos$ & $cos(a)$\\
    $cosh$ & $\frac{(e^{a} + e^{-a})}{2}$ \\
    $exp$ & $e^{a}$\\
    $exp2$ & $2^{a}$\\
    $exp10$ & $10^{a}$\\
    $fmod$ & $a-(floor(\frac{a}{c}) \times c)$\\
    $hypot$ & $\frac{(a^{2} + b^{2})}{2}$\\
    $log$ & $ln(a)$\\
    $log2$ & $\frac{ln(a)}{ln(2)}$\\
    $log10$ & $log(a)$\\
    $pow$ & $a^{c}$\\
    $sin$ & $sin(a)$\\
    $sinh$ & $\frac{(e^{a} - e^{-a})}{2}$\\
    $sqrt$ & $a^{\frac{1}{2}}$ \\
    $tan$ & $tan(a)$\\
    $tanh$ & $\frac{(e^{a} - e^{-a})}{(e^{a} + e^{-a})}$\\
  \end{tabular}
}
  \caption{Formulas for scientific functions in GLIBC. }
  \label{table:formulas}
\end{table}

\subsubsection{Set up the standard of evaluation}

We evaluated our approaches on the functions of 
scientific mathematical library of GLIBC (the GNU C Library), version 2.19. 
This evaluation involves 54 functions 
and we monitor the relative error over $15000$ assignment statements. 
We set the original precision to double and the adjusted precision to 120-bit. 

As described before, in our detecting approach, 
we need to set the threshold of relative error $e0$, 
the threshold of error percentage $p0$
and the threshold for statement percentage $p1$. 
In this situation, an operation, which is executed $m$ times, 
is considered to be a precision-specific operation, 
if it produces relative error greater than $e0$ for over $m \times p0$ input values 
and $m$ is larger than $n \times p1$. 
We evaluate the circumstances when $e0$ is set to $10^{-2}$, $10^{-4}$, $10^{-6}$ and $10^{-8}$. 
For each $e0$, we process the approaches with different values of $p0$. 
We try $p0 = 70\%$ first. If nothing found, then $p0 = 60\%$. Finally $p0 = 50\%$.

\subsection{Evaluation}

\subsubsection{The precision of detecting approach}

We get the evaluation of the precision below in Table~\ref{table:detecting}. 
The precision of detecting approach ranges from $74.51\%$ to $76.58\%$. 
We can figure out the higher the threshold value, the lower the recall rate. 
Besides, it performs better when $e0$ is set to $10^{-6}$ or $10^{-8}$.

\begin{table}[t!]
  \centering
{  \sffamily\small 
  \begin{tabular}{ccccc}
    \textbf{threshold $e0$} & \textbf{$10^{-2}$} & \textbf{$10^{-4}$} & \textbf{$10^{-6}$} & \textbf{$10^{-8}$}\\
    \hline
    \textbf{particular constant} & $59$ & $70$ & $79$ & $79$\\
    \textbf{union} & $6$ & $6$ & $6$ & $6$\\
    \textbf{misjudgement} & $21$ & $26$ & $26$ & $26$\\
    \textbf{total} & $59$ & $70$ & $79$ & $79$\\
    \textbf{precision} & $75.58\%$ & $74.51\%$ & $76.58\%$ & $76.58\%$\\
    \textbf{recall rate} & $73.03\%$ & $85.39\%$ & $95.51\%$ & $95.51\%$\\
  \end{tabular}
}
  \caption{Precision of detecting approach. }
  \label{table:detecting}
\end{table}

\subsubsection{Distribution of precision-specific computation}

In the $54$ tested functions, 
there are $24$ functions are detected to contain precision-specific operations. 
In total, $111$ precision-specific operations are detected. 
After the exclusion of duplicate operations, $43$ precision-specific operations remain unique. 
The distribution in $24$ functions are shown in Table~\ref{table:distribution}

\begin{table}[t!]
  \centering
{  \sffamily\small 
  \begin{tabular}{cc}
    \textbf{Function} & \textbf{Number of }\\
     & \textbf{precision-specific operations} \\
    \hline
    $acosh$ & $1$\\
    $asinh$ & $1$\\
    $atan$ & $1$\\
    $atan2$ & $4$ \\
    $cos$ & $4$\\
    $cosh$ & $2$ \\
    $erf$ & $2$ \\
    $erfc$ & $2$ \\
    $exp$ & $2$\\
    $exp2$ & $1$\\
    $exp10$ & $3$\\
    $j0$ & $10$ \\
    $j1$ & $10$ \\
    $jn$ & $14$ \\
    $lgamma$ & $1$ \\
    $log$ & $1$\\
    $pow$ & $8$\\
    $sin$ & $4$\\
    $sinh$ & $2$\\
    $tan$ & $1$\\
    $tagma$ & $10$ \\
    $y0$ & $7$ \\
    $y1$ & $7$ \\
    $yn$ & $13$ \\
  \end{tabular}
}
  \caption{Formulas for scientific functions in GLIBC. }
  \label{table:distribution}
\end{table}

\subsubsection{Analysis of precision-specific operations}

\paragraph{Particular constant}

This pattern is the most common one. It has a format like:

$x = x1 + constant;$

$y = x – constant;$

A number add a particular constant, and then subtract the particular constant. 
The definition of all the related particular constant is below:

$static\;const\;double\;big\;=\;0x1.8000000000000p45;$

$static\;const\;double\;toint\;=\;0x1.8000000000000p52;$

$const\;static\;mynumber\;three33\;=\;{{0,\;0x42180000}};$

$const\;static\;mynumber\;three51\;=\;{{0,\;0x43380000}};$

$static\;const\;double\;THREEp42\;=\;13194139533312.0;$

$static\;const\;double\;t22\;=\;0x1.8p22;$

$const\;static\;mynumber\;bigu\;=\;{{0xfffffd2c\;,\;0x4297ffff}};$

$const\;static\;mynumber\;bigv\;=\;{{0xfff8016a\;,\;0x4207ffff}};$

The values of these particular constants are listed in Table~\ref{table:constants}

\begin{table}[t!]
  \centering
{  \sffamily\small 
  \begin{tabular}{rl}
    \textbf{Particular constant} & \textbf{Value} \\
    \hline
    $big$ & $1.5 \times 2^{45}$\\
    $toint$ & $1.5 \times 2^{52}$\\
    $three33$ & $1.5 \times 2^{34}$\\
    $three51$ & $1.5 \times 2^{52}$\\
    $THREEp42$ & $1.5 \times 2^{41}$\\
    $t22$ & $1.5 \times 2^{22}$\\
    $bigu$ & $1.5 \times 2^{42} - 724 \times 2^{-10}$\\
    $bigv$ & $1.5 \times 2^{35} - 1 + 362 * 2^{19}$\\
  \end{tabular}
}
  \caption{Values of particular constants. }
  \label{table:constants}
\end{table}

We could notice that $toint$ and $three51$ have a same value. 
The function of these two constants are most obvious. 
They round the floating-point number in double to the nearest integer. 
In the double format, the floating-point number has a sign bit, 
$11$ exponent bits and $52$ fraction bits. 
In order to add, we need the exponents of the two numbers to be the same. 
Then add the two mantissas of the adjusted numbers together. 
Finally adjust the result into standard of floating-point representation. 
For instance, $1101.11$ is the binary representation of $13.75$. We call it $X$. $X$ equals to $1.10111_{(two)} \times 23$.
Under the IEEE754 Standard, it is represented as:

$0,10000000010,1011100000000000000…0000000$

We need the exponents to be same as $1.1_{(two)} \times 252$. 
We call it $Y$. The difference of the exponent is $49$. 
So, add $49$ to $X$’s exponent, and shift the mantissas right by $49$ bit. 
We lose the last two bits $11_{(two)}$ then. 
Because of the rounding, the last $4$ bits, $1101$, turn to $1110$. 
This results in:

$0.0000000000000000…00000000001110_{(two)} \times 2^{52}$

Call this readjusted value, $X’$. 
Next, we add two mantissas of $X’$ and $Y$. The sum is:

$1.10000000000000…0000000000001110_{(two)} \times 2^{52}$

In the IEEE754 Standard, it is represented as:

$0,10000110011,10000000…000000000000001110$

Then, subtract $1.5 \times 2^{52}$

$\;0,10000110011,10000000000…00000000001110$

$- 0,10000110011,1000000000…000000000000000$

$=$

$\;0,10000000010,110000000000000000..000000$

Finally, we get $14$. For the negative number, the situation is similar. 
By shifting the mantissas out, 
the floating-point number in double is rounded to the nearest integer.

Function exp is one of the instances that use particular number $three51$. 
The core source code of the transformed program is below (inserted functions and definitions of variables are omitted):

$temp\_var\_for\_tac\_8\;=\;x\;*\;log2e.x;$

$temp\_var\_for\_tac\_9\;=\;temp\_var\_for\_tac\_8\;+\;three51.x;$

$y\;=\;temp\_var\_for\_tac\_9;$

$temp\_var\_for\_tac\_10\;=\;y\;–\;three51.x;$

When $x = 0.45$, we could trace the relative error as following:

$temp\_var\_for\_tac\_8\_tag5754$

$ORIGINAL:\;		6.49212768400034\;*\;10ˆ-1,\;52/120\;bit$

$SHADOW\;VALUE:\;	6.49212768400034\;*\;10ˆ-1,\;105/120\;bit$

$ABSOLUTE\;ERROR:\;3.82215930115777\;*\;10ˆ-17,\;51/120\;bit$

$RELATIVE\;ERROR:\;5.88737542944107\;*\;10ˆ-17,\;119/120\;bit$

$temp\_var\_for\_tac\_9\_tag5755$

$ORIGINAL:\;		6.75539944105574\;*\;10ˆ15,\;53/120\;bit$

$SHADOW\;VALUE:\;	6.75539944105574\;*\;10ˆ15,\;120/120\;bit$

$ABSOLUTE\;ERROR:\;-3.50787231599966\;*\;10ˆ-1,\;66/120\;bit$

$RELATIVE\;ERROR:\;5.19269415022400\;*\;10ˆ-17,\;118/120\;bit$

$y\_tag5756$

$ORIGINAL:\;		6.75539944105574\;*\;10ˆ15,\;53/120\;bit$

$SHADOW\;VALUE:\;	6.75539944105574\;*\;10ˆ15,\;120/120\;bit$

$ABSOLUTE\;ERROR:\;-3.50787231599966\;*\;10ˆ-1,\;66/120\;bit$

$RELATIVE\;ERROR:\;5.19269415022400\;*\;10ˆ-17,\;118/120\;bit$

$temp\_var\_for\_tac\_10\_tag5757$

$ORIGINAL:\;		1.00000000000000\;*\;10ˆ0,\;1/120\;bit$

$SHADOW\;VALUE:\;	6.49212768400034\;*\;10ˆ-1,\;67/120\;bit$

$ABSOLUTE\;ERROR:\;-3.50787231599966\;*\;10ˆ-1,\;66/120\;bit$

$RELATIVE\;ERROR:\;5.40327067910990\;*\;10ˆ-1,\;120/120\;bit$

After the execution of $y – three51.x$, 
the relative error surges to $5.40327067910990 \times 10ˆ{-1}$ from $10^{-17}$ level. 
We notice the shadow value of $temp\_var\_for\_tac\_104$ equals to 
the shadow value of $temp\_var\_for\_tac\_8$. 
This is correct from the view of real number, 
but disable the rounding function of the particular constant $three51$. 
In other words, high precision keeps the mantissas but is against developers' will. 
Other particular constants have analogous function.

\paragraph{Union}

Two purposes of using union for floating-point numbers are: 
\begin{itemize}[noitemsep]
\item Judge the scale rapidly in condition statement 
\item Modify specific bits concurrently. 
\end{itemize}
For the first situation, the precision-specific operations could not occur, 
since there is no assignment statement. 
For the second case, modification of the integer array may produce large relative error in high precision. 
The difference between the formats of two precisions lead to unforeseen value of the floating-point number.

\subsection{Evaluation of fixing approach}

In order to evaluate our approach, comparison among the results of our approach, 
original precision and higher precision is needed. 
We need four kinds of results of one test program. 
The first one is the set of results of original precision (OP in short). 
The second one is the set of results of higher precision (HP in short), 
which executed with the help of FPdebug. 
The third one is also a set of results of higher precision, 
but with precision-specific operation problem fixed (MP in short). 
The fourth one is a set of standard results obtained 
by aforementioned reliable multi-precision calculation library (S in short). 
For each input, we compute the relative error between OP and S, HP and S, MP and S.

\paragraph{Management of Inaccuracy}

Precondition for the comparison is that OP, HP and MP share the same input. 
This seems easy if we don’t take floating-point inaccuracy into consideration. 
Note that the input of reliable multi-precision calculation library 
is not any floating-point type but a string. 
S is computed under Windows, while OP, HP and MP is computed under Linux. 
This difference is caused by the requirements of 
reliable multi-precision calculation library and FPdebug. 
In the beginning, our experiment simply use a “for” statements to generate inputs, 
and use sprint to convert floating-point type to string. 
Then we realize that this method actually generates three types of inputs! 
Be aware that for statement for a floating-point aggregates inaccuracies along the iterations. 
HP and MP takes more precise inputs, because they calculate the accumulation in higher precision. 
The inputs of S are expected to be equal to those of OP, 
but the conversion by sprintf is not precise enough under Windows. 
To avoid floating-point inaccuracy in the input stage, 
all our inputs are read from files that contains inputs generates by a “for” statement. 
OP, HP and MP read the inputs as floating-point types, while S reads them as strings. 
Floating-point inaccuracies do not only affect our experiment in the input stage, 
but also affect our experiment when we use Excel to calculate 
the average relative error of OP and S, HP and S, MP and S. 
Our results are very sensitive to these inaccuracies, 
so we choose to use reliable multi-precision calculation library 
to process experiment data even if this may take more time.

\paragraph{Comparison of advantage}

We randomly choose $18$ functions to do the evaluation. 
For every input value of each function, we calculate the relative error of MP, HP and LP to S. 
Then we analyze the percentage of the input values that produce less relative error. 
The result is show in Table~\ref{table:performance}

\begin{table*}[t!]
  \centering
{  \sffamily\small 
  \begin{tabular}{ccccccc}
    \textbf{function} & \textbf{$M \geq L$} & \textbf{$M > L$} & \textbf{$M \geq H$} & \textbf{$M > H$} & \textbf{$H \geq L$} & \textbf{$H > L$}\\
    \hline
    \textbf{acos} & $97.00\%$ & $96.90\%$ & $97.00\%$ & $0.20\%$ & $100.0\%$ & $99.90\%$\\
    \textbf{acosh}& $86.91\%$ & $86.81\%$ & $100.00\%$ & $99.00\%$ & $0.60\%$ & $0.50\%$\\
    \textbf{asin} & $96.60\%$ & $96.40\%$ & $96.90\%$ & $0.00\%$ & $99.60\%$ & $99.40\%$\\
    \textbf{asinh} & $83.70\%$ & $83.60\%$ & $100.00\%$ & $95.90\%$ & $1.90\%$ & $1.80\%$\\
    \textbf{atan} & $99.50\%$ & $99.40\%$ & $100.00\%$ & $0.00\%$ & $99.50\%$ & $99.40\%$\\
    \textbf{atan2} & $97.72\%$ & $97.72\%$ & $99.83\%$ & $0.00\%$ & $97.72\%$ & $97.72\%$\\
    \textbf{atanh} & $69.67\%$ & $69.57\%$ & $100.00\%$ & $0.00\%$ & $69.67\%$ & $69.57\%$\\
    \textbf{cos} & $100.00\%$ & $100.00\%$ & $100.00\%$ & $91.91\%$ & $8.09\%$ & $8.09\%$\\
    \textbf{cosh} & $100.00\%$ & $100.00\%$ & $100.00\%$ & $100.00\%$ & $0.00\%$ & $0.00\%$\\
    \textbf{exp} & $100.00\%$ & $100.00\%$ & $100.00\%$ & $100.00\%$ & $0.00\%$ & $0.00\%$\\
    \textbf{exp10} & $100.00\%$ & $100.00\%$ & $100.00\%$ & $100.00\%$ & $0.00\%$ & $0.00\%$\\
    \textbf{exp2} & $100.00\%$ & $99.90\%$ & $100.00\%$ & $99.90\%$ & $0.10\%$ & $0.00\%$\\
    \textbf{log} & $99.98\%$ & $99.98\%$ & $100.00\%$ & $89.70\%$ & $10.28\%$ & $10.28\%$\\
    \textbf{log10} & $88.90\%$ & $88.90\%$ & $100.00\%$ & $0.00\%$ & $88.90\%$ & $88.90\%$\\
    \textbf{log2} & $94.90\%$ & $94.90\%$ & $100.00\%$ & $0.00\%$ & $94.90\%$ & $94.90\%$\\
    \textbf{sin}  & $100.00\%$ & $100.00\%$ & $100.00\%$ & $92.04\%$ & $7.96\%$ & $7.96\%$\\
    \textbf{sinh} & $83.60\%$ & $83.50\%$ & $100.00\%$ & $56.00\%$ & $27.60\%$ & $27.50\%$\\
    \textbf{tan} & $68.81\%$ & $68.81\%$ & $100.00\%$ & $74.94\%$ & $18.22\%$ & $18.22\%$\\
    \textbf{tanh} & $99.50\%$ & $21.38\%$ & $100.00\%$ & $0.00\%$ & $99.50\%$ & $21.38\%$\\
  \end{tabular}
}
  \caption{Percentage of the input values that perform better. }
  \label{table:performance}
\end{table*}

The results shown in the table is astonishing, 
that there may be no advantage to perform all the operations in high precision. 
Except function asin, the percentages of the input values with 
better performance in high precision than original precision are all less than $30\%$. 
Our common sense, that rising precision results in better accuracy, may not be a fact.
Furthermore, our fixing method has significant advantage. 
In nearly all the functions, the percentages of the input values 
with better performance in our method than high precision or original precision are over $90\%$. 
In addition, we could see that precision-specific operation influences a lot to the final result. 
Our method are mainly processed in high precision, 
with decreasing precision for precision-specific operations. 
This simple change make the performance much better. 
For example, in function $exp$, no input are better in HP than OP. 
After applying the fixing approach, $100.00\%$ perform better.

\paragraph{Comparison of average}

We calculate the average of relative errors for each function in MP, HP and LP. 
The result is shown in Table~\ref{table:average}

\begin{table*}[t!]
  \centering
{  \sffamily\small 
  \begin{tabular}{cccc}
    \textbf{function} & \textbf{$AVG MP$} & \textbf{$AVG LP$} & \textbf{$AVG HP$}\\
    \hline
    \textbf{acos} & $5.899867786E-11$ & $4.104786918E-17$ & $9.726528892E-20$\\
    \textbf{acosh} & $1.032899398E-17$ & $4.358485412E-17$ & $5.512172898E+04$\\
    \textbf{asin} & $4.031172822E-11$ & $3.989935789E-17$ & $2.089788145E-19$\\
    \textbf{asinh} & $1.300417912E-17$ & $4.356108197E-17$ & $7.535137885E+03$\\
    \textbf{atan} & $4.023136878E-19$ & $3.587342031E-17$ & $4.023136878E-19$\\
    \textbf{atan2} & $1.400411860E-18$ & $4.084347768E-17$ & $1.394693261E-18$\\
    \textbf{atanh} & $3.329345125E-17$ & $4.890462024E-17$ & $3.329345125E-17$\\
    \textbf{cos} & $2.328553719E-21$ & $3.710995528E-17$ & $2.595204983E-03$\\
    \textbf{cosh} & $4.290083004E-22$ & $4.415740936E-17$ & $9.541710311E-07$\\
    \textbf{exp} & $1.908302538E-25$ & $4.126635952E-17$ & $9.540155830E-07$\\
    \textbf{exp10} & $8.484205894E-23$ & $6.977866446E-17$ & $1.262973185E-06$\\
    \textbf{exp2} & $6.029968486E-20$ & $4.116251294E-17$ & $3.249127529E-04$\\
    \textbf{log} & $8.619452474E-21$ & $4.082200343E-17$ & $1.091593070E+01$\\
    \textbf{log10} & $8.117175153E-18$ & $4.265676679E-17$ & $8.117175153E-18$\\
    \textbf{log2} & $6.350986821E-18$ & $4.364585535E-17$ & $6.350986821E-18$\\
    \textbf{sin} & $7.735410706E-21$ & $3.565606085E-17$ & $4.273276871E-02$\\
    \textbf{sinh} & $1.767186452E-17$ & $4.653248592E-17$ & $5.331589632E-07$\\
    \textbf{tan} & $1.713519956E-17$ & $3.893237360E-17$ & $1.433518223E+13$\\
    \textbf{tanh} & $1.689024716E-19$ & $5.969240750E-18$ & $1.689024716E-19$\\
  \end{tabular}
}
  \caption{Average relative error to standard value. }
  \label{table:average}
\end{table*}

The comparison of average relative errors shows for most functions, 
our fixing method could reduce the relative error significantly. 
Because of the precision-specific operation problem, 
rising all operations into high precision is not a good choice.

%% file: Conclusion.tex
\section{Conclusion}

In this paper, we propose a light-weight detecting approach for precision-specific operation 
and a very efficient fixing approach based on the adjustment of precision. 
We evaluate our method on the basic scientific mathematical library of GLIBC. 
The result proves our method are more accurate than original precision and high precision. 
Contrary to the common sense, we discover that simply rising precision for all operations, 
without fixing the precision-specific operation problem, may lead to larger relative error.